\documentclass[12pt]{article}
\usepackage[matrix,arrow]{xy}
\usepackage{amssymb,amsfonts,amsmath,footnote}

\newtheorem{thm}{Th\'eor\`eme}[section]

\newtheorem{lem}[thm]{Lemme}

\newtheorem{prop}[thm]{Proposition}

\newtheorem{Def}[thm]{D\'efinition}

\newtheorem{conj}[thm]{Conjecture}

\newtheorem{exem}[thm]{Exemple}

\numberwithin{equation}{section}

\def\NN{\mathbb N}
\def\CCC{\mathbb C}
\def\RRR{\mathbb R}

\def\EE{\mathbb E}

\def\oge{\leavevmode\raise
.3ex\hbox{$\scriptscriptstyle\langle\!\langle\,$}}
\def\feg{\leavevmode\raise
.3ex\hbox{$\scriptscriptstyle\,\rangle\!\rangle$}}

\date{S\'eminaire Bourbaki, Juin 2014; \\
66\`eme ann\'ee, 2013-2014; $n^o$ 1088}
\title{LE PROBL\`EME DE KADISON-SINGER (d'apr\`es A. Marcus, D. Spielman et N. Srivastava)}
\author{Alain VALETTE}

\begin{document}
\maketitle

{\it The results that we have obtained leave the question of uniqueness of extension of the singular pure states of ${\mathcal A}_d$ open. We incline to view that such extension is non-unique (...)} (R.V. Kadison and I.M. Singer, \cite{KS}, p.~397). 

\bigskip
\section{Introduction}

\bigskip

\subsection{\'Enonc\'e du probl\`eme}

Si $\mathcal{H}$ est un espace de Hilbert, une {\it $C^*$-alg\`ebre} sur $\mathcal{H}$ est une $*$-sous-alg\`ebre, ferm\'ee pour la norme-op\'erateur, de l'alg\`ebre $\mathcal{B}(\mathcal{H})$ des op\'erateurs lin\'eaires born\'es sur $\mathcal{H}$. Si $A$ est une $C^*$-alg\`ebre unitale sur $\mathcal{H}$, un {\it \'etat} sur $A$ est une fonctionnelle lin\'eaire continue $\varphi$ sur $A$, positive (c-\`a-d. $\varphi(x^*x)\geq 0$ pour tout $x\in A$), et normalis\'ee ($\varphi(1)=1$). L'espace $S(A)$ des \'etats sur $A$ est une partie convexe faible-$*$ compacte du dual de $A$; un \'etat est {\it pur} si c'est un point extr\^eme de $S(A)$. 

\begin{exem}\label{vectoriel} Avec $A= \mathcal{B}(\mathcal{H})$: si $\xi\in\mathcal{H}$ est un vecteur de norme 1, l'\'etat $\varphi(T)=\langle T(\xi)|\xi\rangle$ est pur; ces \'etats s'appellent {\it \'etats vectoriels}.
\end{exem}

\begin{exem} Si $A$ est commutative, les \'etats purs sont exactement les caract\`eres de $A$. En effet, la transform\'ee de Gelfand identifie $A$ \`a $C(X)$, la $C^*$-alg\`ebre des fonctions continues sur un espace compact $X$; par le th\'eor\`eme de repr\'esentation de Riesz, $S(C(X))$ s'identifie \`a l'espace des mesures de probabilit\'es sur $X$; une telle mesure est extr\'emale si et seulement si c'est une masse ponctuelle, donc l'\'etat correspondant est l'\'evaluation en un point de $X$.
\end{exem}

Pour $A$ une $C^*$-alg\`ebre unitale de $\mathcal{B}(\mathcal{H})$, une application du th\'eor\`eme de Hahn-Banach montre que tout \'etat $\varphi$ de $A$ s'\'etend en un \'etat de $\mathcal{B}(\mathcal{H})$ (voir \cite{Dix}, lemme 2.10.1). L'ensemble des extensions de $\varphi$ est un convexe $K_\varphi$ faible-$*$ compact de $\mathcal{B}(\mathcal{H})^*$. Si $\varphi$ est pur, les points extr\^emes de $K_\varphi$ sont des \'etats purs de $\mathcal{B}(\mathcal{H})$ (en effet, dans ce cas $K_\varphi$ est une face de $S(\mathcal{B}(\mathcal{H}))$). Par le th\'eor\`eme de Krein-Milman, un \'etat pur sur $A$ s'\'etend en un unique \'etat pur de $\mathcal{B}(\mathcal{H})$ si et seulement s'il s'\'etend en un unique \'etat de $B(H)$.

\begin{exem} $\mathcal{B}(\mathcal{H}))$ admet des \'etats purs qui ne sont pas vectoriels. En effet, consid\'erons $C[0,1]$ comme une $C^*$-alg\`ebre agissant par multiplication sur $\mathcal{H}$$=L^2[0,1]$. Un \'etat pur de $\mathcal{B}(\mathcal{H})$ qui \'etend un \'etat pur de $C[0,1]$ ne peut \^etre vectoriel, car les fonctions de $C[0,1]$ n'ont pas de vecteur propre commun sur $\mathcal{H}$.
\end{exem}

\begin{exem} Il est facile de voir qu'une extension d'un \'etat non pur n'est pas unique en g\'eneral. Prenons par exemple l'alg\`ebre $A\simeq\CCC^2$ des matrices diagonales dans l'alg\`ebre $M_2(\CCC)$ des matrices 2-fois-2 complexes. Consid\'erons les matrices positives $B_1=\left(\begin{array}{cc}1/2 & 0 \\0 & 1/2\end{array}\right)$ et $\left(\begin{array}{cc}1/2 & i/2 \\-i/2 & 1/2\end{array}\right)$. Les \'etats $S\mapsto Tr(SB_1)$ et $S\mapsto Tr(SB_2)$ de $M_2(\CCC)$ se restreignent tous deux en l'\'etat non pur $(a_1,a_2)\mapsto\frac{a_1+a_2}{2}$ de $A$.
\end{exem}

En 1959, R.V. Kadison et I.M. Singer \cite{KS} \'etudient la propri\'et\'e d'extension unique des \'etats purs sur les $C^*$-alg\`ebres ab\'eliennes maximales (MASA) de $\mathcal{B}(\mathcal{H})$. Si $\mathcal{H}$=$L^2[0,1]$ (avec la mesure de Lebesgue) et $A=L^\infty[0,1]$ agissant par multiplication, ils montrent que $A$ n'a pas la propri\'et\'e d'extension unique des \'etats purs (\cite{KS}, Theorem 2). Pour $\mathcal{H}$=$\ell^2(\NN)$ et l'alg\`ebre $D=\ell^\infty(\NN)$ des op\'erateurs diagonaux, ils posent la question: 

\medskip
{\bf Probl\`eme de Kadison-Singer: $D$ a-t-elle la propri\'et\'e d'extension unique des \'etats purs?}

\begin{exem} Pour $k\in\NN$, l'\'etat pur $\varphi_k$ de $D$, d\'efini par $\varphi_k((a_n)_{n\in\NN})=a_k$, admet comme unique extension \`a $\mathcal{B}$$(\ell^2(\NN))$ l'\'etat pur $T\mapsto\langle T{\bf e}_k|{\bf e}_k\rangle$, o\`u ${\bf e}_k$ est le $k$-i\`eme vecteur de la base canonique de $\ell^2(\NN)$. En effet, soit $\psi$ un \'etat pur de $\mathcal{B}$$(\ell^2(\NN))$ qui \'etend $\varphi_k$, et $\pi_\psi$ la repr\'esentation GNS associ\'ee. L'\'etat $\psi$ ne s'annule pas sur l'id\'eal $\mathcal{K}$ des op\'erateurs compacts (puisqu'il prend la valeur 1 sur le projecteur de rang 1 associ\'e \`a ${\bf e}_k$). Ainsi la restriction de $\pi_\psi$ \`a $\mathcal{K}$ est une repr\'esentation irr\'eductible de $\mathcal{K}$, elle est donc \'equivalente \`a la repr\'esentation standard de $\mathcal{K}$, donc $\psi$ est un \'etat vectoriel; c'est l'\'etat associ\'e \`a~${\bf e}_k$, comme on le voit facilement\footnote{Notons que $D\simeq C(\beta\NN)$, o\`u $\beta\NN$ est le compactifi\'e de Stone-$\check{C}$ech de $\NN$, c-\`a-d. l'ensemble des ultrafiltres de $\NN$. Pour attaquer le probl\`eme de Kadison-Singer, il suffirait donc de consid\'erer les \'etats purs de $D$ correspondant aux ultrafiltres libres de $\NN$. Cette approche a \'et\'e consid\'er\'ee (voir par exemple \cite{AndTAMS}), mais n'intervient pas dans la r\'ecente solution du probl\`eme de Kadison-Singer.}.
\end{exem}

Il y a une fa\c con canonique d'\'etendre un \'etat de $D$ en un \'etat de $\mathcal{B}(\ell^2(\NN))$, gr\^ace \`a l'esp\'erance conditionnelle $E:\mathcal{B}(\ell^2(\NN))$$\rightarrow D: T\mapsto diag(T)$ (o\`u $diag(T)$ est la diagonale de $T$ dans la base canonique): si $\varphi$ est un \'etat de $D$, alors $\varphi\circ E$ est un \'etat de $\mathcal{B}(\ell^2(\NN))$ qui \'etend $\varphi$ (si $\varphi$ est pur, alors $\varphi\circ E$ est un \'etat pur de $\mathcal{B}(\ell^2(\NN))$, voir \cite{AndJFA}, Theorem 1). Le probl\`eme de Kadison-Singer se reformule donc ainsi~: si $\varphi$ est un \'etat pur de $D$, l'\'etat $\varphi\circ E$ est-il la seule extension en un \'etat de $\mathcal{B}(\ell^2(\NN))$? Quoiqu'ils aient \'etabli l'unicit\'e de~$E$ comme esp\'erance conditionnelle $\mathcal{B}(\ell^2(\NN))$$\rightarrow D$ (voir \cite[Thm.1]{KS}), Kadison et Singer tendaient \`a penser que la r\'eponse devait \^etre n\'egative, comme l'indique la phrase en exergue\footnote{Le probl\`eme d'extension unique des \'etats purs peut se reformuler pour d'autres MASAs dans d'autres alg\`ebres de von Neumann. En particulier, dans un article r\'ecent \cite{Pop}, Popa montre que le probl\`eme de Kadison-Singer est \'equivalent \`a la propri\'et\'e d'extension unique des \'etats purs pour l'ultraproduit des alg\`ebres diagonales dans l'ultraproduit des $M_n(\CCC)$, qui est un facteur de type $II_1$; il montre \'egalement, ind\'ependamment des r\'esultats de Marcus-Spielman-Srivastava, que si $\omega$ est un ultrafiltre libre sur $\NN$ et si $A$ est une MASA singuli\`ere dans un facteur $M$ de type $II_1$, alors l'ultraproduit $A^\omega$ a la propri\'et\'e d'extension unique des \'etats purs dans l'ultraproduit $M^\omega$.}. 

\subsection{Motivations}

Dans \cite{KS}, Kadison et Singer \'ecrivent qu'ils ont appris le probl\`eme de I. Kaplansky et I. Segal. Cependant, lors d'expos\'es plus r\'ecents (voir par exemple \cite{Jor}), Kadison a affirm\'e que le probl\`eme trouve sa source dans certains passages du livre de P.A.M.~Dirac \cite{Dir} sur les fondements de la m\'ecanique quantique. La question trait\'ee par Dirac est: comment d\'eterminer les probabilit\'es {\oge de base\feg} associ\'ees aux \'etats quantiques d'un syst\`eme~? \`A la section 18 de \cite{Dir}, Dirac sugg\`ere de commencer par consid\'erer un ensemble complet d'observables commutant deux \`a deux\footnote{Il s'agit donc d'observables compatibles, telles que la mesure de l'un n'influence pas la mesure de l'autre.}; un tel ensemble forme la partie auto-adjointe d'une MASA. Ensuite, on sp\'ecifie les distributions de probabilit\'e associ\'ees aux observables de cet ensemble dans un \'etat quantique donn\'e: cela revient \`a fixer un \'etat $\varphi$ sur la MASA\footnote{La probabilit\'e que, dans cet \'etat, une observable $A$ prenne des valeurs entre $a$ et $b$ est alors donn\'ee par $\varphi(E_{[a,b])})$, o\`u $E_{[a,b]}$ est le projecteur spectral de $A$ associ\'e \`a l'intervalle $[a,b]$.}. Enfin, on veut g\'en\'eraliser cette information au syst\`eme entier, c-\`a-d. \'etendre cette distribution de probabilit\'e \`a toutes les autres observables, m\^eme non compatibles avec celles de la MASA. Il est naturel de se demander si cette extension est unique, ce qui am\`ene au probl\`eme de Kadison-Singer\footnote{Citons la p.~75 de \cite{Dir}: {\oge {\it The representation is then completely determined except for the arbitrary phase factors. For most purposes the abritrary phase factors are unimportant and trivial, so that we may count the representation as being completely determined by the observables that are diagonal in it.}\feg} Cela peut laisser penser que Dirac \'etait convaincu de l'unicit\'e de l'extension d'un \'etat. Nous renvoyons \`a \cite{Brian} pour davantage de d\'etails quant \`a l'interpr\'etation physique.}.

\subsection{Les principaux r\'esultats}

Le probl\`eme de Kadison-Singer a \'et\'e r\'esolu -- par l'affirmative ! -- en juin 2013 par A.~Marcus, D. Spielman et N. Srivastava \cite{MSS-KS}, trois scientifiques actifs \`a l'interface des math\'ematiques et de l'informatique th\'eorique:

\begin{thm}\label{KS-main} Tout \'etat pur sur $D$ admet une unique extension en un \'etat de ${\mathcal B}(\ell^2(\NN))$ (n\'ecessairement pur). 
\end{thm}

Ce r\'esultat n'a pas \'et\'e d\'emontr\'e directement, mais via une traduction en alg\`ebre lin\'eaire  due \`a J. Anderson \cite{AndTAMS}, C. Akemann et J. Anderson \cite{AkA}, et N. Weaver \cite{Wea}. Cette traduction est donn\'ee \`a la section \ref{finKS}. Les principaux r\'esultats de A. Marcus, D. Spielman et N. Srivastava \cite[Thm.1.2]{MSS-KS} sont les suivants (o\`u on note $p_A(z)=\det(z.1_m-A)$ le polyn\^ome caract\'eristique de $A$, pour une matrice $A$ de taille $m$).

\begin{thm}\label{MSS-expect} Soient $A_1,...,A_d$ des variables al\'eatoires ind\'ependantes \`a valeurs dans les matrices semi-d\'efinies positives de rang $1$ de $M_m(\CCC)$; posons $A=\sum_{i=1}^d A_i$. On suppose que $\EE A=1_m$ et que $\EE \|A_i\|\leq\varepsilon$ pour $i=1,...,d$. Alors la plus grande racine r\'eelle de $\EE p_A$ est au plus $(1+\sqrt{\varepsilon})^2$.
\end{thm}

\begin{thm}\label{MSS-posproba} Soient $A_1,...,A_d$ des variables al\'eatoires ind\'ependantes \`a valeurs dans les matrices semi-d\'efinies positives de rang $1$ de $M_m(\CCC)$; posons $A=\sum_{i=1}^d A_i$. On suppose que les $A_i$ prennent un nombre fini de valeurs. Alors $\|A\|$ est inf\'erieure ou \'egale \`a la plus grande racine r\'eelle de $\EE p_A$, pour au moins une r\'ealisation des $A_i$.
\end{thm} 

Comme $\|A\|$ est la plus grande racine positive de $p_A$, le th\'eor\`eme \ref{MSS-posproba} se reformule: pour au moins une r\'ealisation des $A_i$, la plus grande racine de $p_A$ est inf\'erieure ou \'egale \`a la plus grande racine de $\EE p_A$. On peut donc penser \`a ce r\'esultat comme une version {\oge non lin\'eaire\feg} d'un principe de probabilit\'e \'el\'ementaire: si $X$ est une variable al\'eatoire r\'eelle prenant un nombre fini de valeurs, alors $X\leq \EE X$ pour au moins une r\'ealisation de X.

\subsection{\'Enonc\'es \'equivalents}

Une s\'erie de travaux ont montr\'e, depuis la fin des ann\'ees 1970, que le probl\`eme de Kadison-Singer est \'equivalent \`a une s\'erie d'\'enonc\'es ouverts en alg\`ebre lin\'eaire, en analyse harmonique, en th\'eorie des op\'erateurs, et en analyse du signal - \'enonc\'es qui se retrouvent donc d\'emontr\'es du m\^eme coup! Citons:
\begin{itemize}
\item la conjecture $(KS_r)$ de N. Weaver \cite{Wea}; voir les remarques pr\'ec\'edant le th\'eor\`eme \ref{Weaver} ci-dessous;
\item la conjecture de pavage de J. Anderson \cite{AndTAMS}; voir le th\'eor\`eme \ref{Anderson} ci-dessous;
\item la conjecture de Bourgain-Tzafriri sur les matrices dont toutes les colonnes sont de norme 1 - voir la section \ref{BTF}.
\item la conjecture de Feichtinger en th\'eorie des frames\footnote{Nous utilisons cet anglicisme avec la b\'en\'ediction de S. Mallat \cite{Mal}.} -- voir \'egalement la section \ref{BTF}.
\end{itemize}
Les \'equivalences entre ces \'enonc\'es  sont bien expliqu\'ees dans les travaux de P. Casazza et de ses collaborateurs (voir par exemple \cite{CFTW}). 

\subsection{Organisation de l'article}

Les th\'eor\`emes \ref{MSS-expect} et \ref{MSS-posproba} seront d\'emontr\'es \`a la section \ref{carmixte}. La principale id\'ee est de consid\'erer $p_A$ comme polyn\^ome caract\'eristique mixte de $A_1,..., A_d$ (voir la d\'efinition \ref{polmixte} ci-dessous), qui est une sp\'ecialisation d'un polyn\^ome \`a $d+1$ variables obtenu \`a partir de $\det(z.1_m +\sum_{i=1}^d z_iA_i)$. On est ainsi conduit \`a \'etudier une classe de polyn\^omes \`a plusieurs variables, les polyn\^omes {\it r\'eels stables}, des polyn\^omes \`a coefficients r\'eels d\'efinis par une condition sur les z\'eros complexes (cf. D\'efinition \ref{polresta}), qui sont \'etudi\'es via un m\'elange de m\'ethodes r\'eelles et complexes; cette \'etude est faite \`a la section \ref{reelstable}. La preuve originale de \cite{MSS-KS} utilisait des propri\'et\'es d'entrelacements de z\'eros pour des familles de polyn\^omes. Dans \cite{Tao}, T. Tao remplace les entrelacements de z\'eros par des arguments bas\'es sur le fait que certaines combinaisons convexes de polyn\^omes ont tous leurs z\'eros r\'eels\footnote{On sait que les deux approches sont \'equivalentes, voir \cite[Thm. 2']{Fell}; \cite[Thm.4]{Sr2}.}. Nous nous inspirons tr\`es largement de l'approche de Tao. La solution du probl\`eme de Kadison-Singer est pr\'esent\'ee \`a la section \ref{finKS}; les applications aux conjectures de Feichtinger et Bourgain-Tzafriri sont dans la section \ref{BTF}.

 Marcus, Spielman et Srivastava ne se sont pas content\'es de r\'esoudre le probl\`eme de Kadison-Singer! Deux mois plus t\^ot, ils avaient r\'esolu par l'affirmative un probl\`eme \`a peine moins c\'el\`ebre: celui de l'existence de familles infinies de graphes de Ramanujan de degr\'e arbitraire (voir \cite{MSS-Ra}). Comme il s'agit d'une belle application du th\'eor\`eme \ref{MSS-posproba}, nous donnons leur preuve \`a la section \ref{Raman}.

\section{Notations}
\subsection{Espaces de Hilbert}
On notera $\ell^2=\ell^2(\NN)$ et $\ell^2_m=\ell^2(\{1,...,m\})$; on note ${\bf e}_k$ le $k$-i\`eme vecteur de la base canonique; pour $m<n$, on voit $\ell^2_m$ comme le sous-espace de $\ell^2_n$ engendr\'e par ${\bf e_1,...,e_m}$.

\noindent
Produit scalaire: $\langle{\bf x|y}\rangle=\sum_n x_n\overline{y_n}$.

\noindent
Normes: $\|{\bf x}\|_2=(\sum_n |x_n|^2)^{1/2}$; $\|{\bf x}\|_\infty=\sup_n |x_n|$.

\subsection{Op\'erateurs et matrices} 
On note $\mathcal{B}$$(\ell^2)$ l'alg\`ebre des op\'erateurs lin\'eaires born\'es; on les voit comme matrices de taille infinie par rapport \`a la base $({\bf e}_k)_k$; pour $T\in\mathcal{B}$$(\ell^2)$, on note $diag(T)$ la suite des valeurs diagonales de $T$. La {\it norme-op\'erateur} de $T$ est:
$$\|T\|=\sup_{\|{\bf x}\|_2\leq 1} \|T({\bf x})\|_2.$$

\noindent
Un op\'erateur $P\in\mathcal{B}$$(\ell^2)$ est un {\it projecteur} si $P=P^*=P^2$. Un op\'erateur $T$ est {\it semi-d\'efini positif} si $\langle T({\bf x)|x}\rangle\geq 0$ pour tout ${\bf x}\in\ell^2$; de fa\c con \'equivalente, il existe un op\'erateur $S$ tel que $T=S^*S$. Pour $S,T$ des op\'erateurs auto-adjoints, on note $S\leq T$ si $T-S$ est semi-d\'efini positif.

\noindent
On note $M_m(\CCC)$ l'alg\`ebre $\mathcal{B}$$(\ell^2_m)$ des matrices complexes de taille $m$; la matrice-identit\'e se note $1_m$. Le polyn\^ome caract\'eristique de la matrice $T\in M_m(\CCC)$ est $p_T(z)=\det(z.1_m-T)$. 

\subsection{Fonctions et polyn\^omes}
Si $p(z)$ est un polyn\^ome en une variable, $ZM(p)$ d\'esigne le plus grand z\'ero r\'eel de $p$ (s'il existe).
Pour les fonctions $f(z_1,...,z_d)$ sur $\CCC^d$, on note $\partial_i=\frac{\partial}{\partial z_i}$ l'op\'erateur de d\'eriv\'ee partielle par rapport \`a la $i$-i\`eme variable. On note $\Phi^i_f=\frac{\partial_i f}{f}$ la d\'eriv\'ee logarithmique de $f$ par rapport \`a la $i$-i\`eme variable.

\section{Polyn\^omes r\'eels stables}\label{reelstable}

\subsection{D\'efinitions et premi\`eres propri\'et\'es}

Notons $\mathbb{H}=\{z\in\CCC:Im(z)>0\}$ le demi-plan sup\'erieur de $\CCC$.

\begin{Def}\label{polresta} Un polyn\^ome en $d$ variables $z_1,...,z_d$ est r\'eel stable si ses coefficients sont r\'eels et qu'il n'admet aucun z\'ero dans l'ouvert $\mathbb{H}^d$ de $\CCC^d$.
\end{Def}

En particulier, un polyn\^ome \`a une variable est r\'eel stable, si et seulement si ses coefficients sont r\'eels et tous ses z\'eros sont r\'eels. Voici la principale source d'exemples de polyn\^omes r\'eels stables. 

\begin{lem}\label{mainex} Soient $A_1,...,A_d$ des matrices semi-d\'efinies positives dans $M_m(\CCC)$. Le polyn\^ome $q(z,z_1,...,z_d)=\det(z.1_m+\sum_{i=1}^d z_iA_i)$ est r\'eel stable.
\end{lem}

{\bf Preuve:} En utilisant $A_i=A_i^*$, on voit facilement que $$q(\overline{z},\overline{z_1},...,\overline{z_d})=\overline{q(z,z_1,...,z_d)},$$ donc $q$ est \`a coefficients r\'eels. Supposons par l'absurde que $q$ s'annule en $(z,z_1,...,z_d)$, avec $Im(z)>0$ et $Im(z_i)>0$ pour $i=1,...,d$. On trouve donc un vecteur $\xi\in\ell^2_m$ non nul, dans le noyau de l'op\'erateur $z.1_m+\sum_{i=1}^d z_iA_i$. En particulier $0=z\|\xi\|^2_2+\sum_{i=1}^d z_i\langle A_i(\xi)|\xi\rangle$, et les coefficients de cette combinaison lin\'eaire de $z,z_1,...,z_d$ sont positifs et non tous nuls. En prenant les parties imaginaires, on a une contradiction.
\hfill$\square$

\medskip
Une propri\'et\'e importante de la classe des polyn\^omes r\'eels stables est d'\^etre pr\'eserv\'ee par sp\'ecialisation \`a des valeurs r\'eelles des variables.

\begin{prop}[\cite{Wag}, lemma 2.4(d)] 
\label{special}
Soit $p(z_1,...,z_d)$ un polyn\^ome r\'eel stable, avec $d>1$. Soit $t\in\RRR$. Le polyn\^ome \`a $d-1$ variables $p(z_1,...,z_{d-1},t)$ est r\'eel stable, s'il n'est pas identiquement nul.
\end{prop}

{\bf Preuve:} Il est clair que $q(z_1,...,z_{d-1})=p(z_1,...,z_{d-1},t)$ est \`a coefficients r\'eels. Comme la suite de polyn\^omes $(p(z_1,...,z_{d-1},t+\frac{i}{n}))_{n\geq 1}$ converge vers $q$ uniform\'ement sur tout compact de $\mathbb{H}^{d-1}$ et n'a pas de z\'ero dans $\mathbb{H}^{d-1}$, le th\'eor\`eme de Hurwitz en analyse complexe garantit que $q$ est sans z\'ero dans $\mathbb{H}^{d-1}$, \`a moins d'\^etre identiquement nul. \hfill$\square$

\medskip
La classe des polyn\^omes r\'eels stables est aussi invariante par certains op\'erateurs diff\'erentiels:

\begin{prop}[\cite{MSS-KS}, Cor. 3.7 ; \cite{Wag}, lemma 2.4(f)] 
\label{deriv}
Soit $p(z_1,...,z_d)$ un polyn\^ome r\'eel stable. Pour $t\in\RRR$, le polyn\^ome $(1+t\partial_d)p(z_1,...,z_d)$ est r\'eel stable.
\end{prop}

{\bf Preuve:} On peut bien s\^ur supposer $t\neq 0$. Supposons par l'absurde que $(1+t\partial_d)p$ ait un z\'ero $(z_1,...,z_d)$ dans $\mathbb{H}^d$. Le polyn\^ome $q(z)=p(z_1,...,z_{d-1},z)$ est sans z\'ero dans $\mathbb{H}$, en particulier $q(z_d)\neq 0$. Factorisons $q$ en facteurs du 1er degr\'e: $q(z)=c\prod_{i=1}^n(z-w_i)$, o\`u les z\'eros $w_i$ satisfont $Im(w_i)\leq 0$. On a alors 
$$0=(1+t\partial_d)p(z_1,...,z_d)=(q+tq')(z_d)=q(z_d)(1+t(\frac{q'}{q})(z_d))\,;$$
donc en prenant la d\'eriv\'ee logarithmique de $q$:
$$0=1+t(\frac{q'}{q})(z_d)=1+t \sum_{i=1}^n \frac{1}{z_d-w_i}=1+t\sum_{i=1}^n \frac{\overline{z_d-w_i}}{|z_d-w_i|^2}.$$
En prenant les parties imaginaires:
$$0=t\sum_{i=1}^n \frac{Im(w_i)-Im(z_d)}{|z_d-w_i|^2}.$$
Mais ceci est une contradiction car $t\neq 0$ et $Im(w_i)<Im(z_d)$ pour $i=1,...,d$.
\hfill$\square$

\subsection{Propri\'et\'es de convexit\'e}

\begin{lem} [\cite{Tao}, lemma 14] 
\label{lem14}
Soient $p(z)$ et $q(z)$ deux polyn\^omes en une variable, de m\^eme degr\'e, avec coefficient de plus haut degr\'e \'egal \`a $1$. On suppose que, pour tout $t\in[0,1]$, le polyn\^ome $(1-t)p+tq$ est r\'eel stable. Alors $ZM((1-t)p+tq)$ est dans l'intervalle $[ZM(p),ZM(q)]$.
\end{lem}

{\bf Preuve:} On peut bien s\^ur supposer $0<t<1$ et $ZM(p)\leq ZM(q)$.
\begin{itemize}
\item Commen\c cons par montrer que $ZM((1-t)p+tq)\leq ZM(q)$. Si $x>ZM(q)$, on a simultan\'ement $q(x)>0$ et $p(x)>0$, donc $((1-t)p+tq)(x)>0$, et donc $x>ZM((1-t)p+tq)$.
\item Pour montrer $ZM(p)\leq ZM((1-t)p+tq)$, il suffit de montrer que $q(ZM(p))\leq 0$ (alors $((1-t)p+tq)(ZM(p))\leq 0$, et le r\'esultat d\'ecoule du th\'eor\`eme des valeurs interm\'ediaires).

Supposons par l'absurde que $q(ZM(p))>0$, donc aussi $q(x)>0$ pour $ZM(p)<x<ZM(p)+\varepsilon$, avec $\varepsilon>0$ assez petit.

On va utiliser le caract\`ere r\'eel stable comme suit: la fonction $N:[0,1]\rightarrow\NN$ qui compte le nombre de z\'eros (avec multiplicit\'es) de $(1-t)p+tq$ dans l'intervalle $[ZM(p)+\varepsilon,ZM(q)]$ est constante. En effet, comme $(1-t)p+tq$ est strictement positif \`a gauche de $ZM(p)+\varepsilon$ et \`a droite de $ZM(q)$, les z\'eros ne peuvent pas {\oge s'\'echapper\feg} \`a gauche ou \`a droite. Donc la seule fa\c con pour $N$ de varier est qu'il existe une valeur de $t$ o\`u deux z\'eros r\'eels de $(1-t)p+tq$ se confondent et bifurquent en des z\'eros complexes conjugu\'es (ou vice-versa), ce qui est exclu par le caract\`ere r\'eel stable de $(1-t)p+tq$. Ceci fournit la contradiction d\'esir\'ee, puisque $N(0)=0$ ($p$~est strictement positif \`a droite de $ZM(p))$ et $N(1)\geq 2$ (comme \mbox{$q(ZM(p)+\varepsilon)>0$,} le polyn\^ome $q$ a au moins 2 z\'eros sur l'intervalle consid\'er\'e.
\hfill$\square$
\end{itemize}

\medskip
Pour ${\bf x}=(x_1,...,x_d)\in\RRR^d$, on d\'efinit l'{\it orthant} $\{{\bf y\geq x}\}$ comme:
$$\{{\bf y\geq x}\}=\{{\bf y}=(y_1,...,y_d)\in\RRR^d: y_i\geq x_i\;\mbox{pour tout}\;i=1,...,d\}.$$
Rappelons que, pour une fonction $f$ de $d$ variables, on note $\Phi^j_f=\frac{\partial_jf}{f}$ la d\'eriv\'ee logarithmique de $f$ par rapport \`a la $j$-i\`eme variable ($j=1,...,d$).

\begin{lem}[\cite{MSS-KS}, lemma 5.7; \cite{Tao}, lemma 17]
\label{lem16+17}
 Soit $p(z_1,...,z_d)$ un polyn\^ome r\'eel stable; soit ${\bf x}=(x_1,...,x_d)\in\RRR^d$; si $p$ est sans z\'ero dans l'orthant $\{{\bf y\geq x}\}$; alors, pour tous $k\in\NN ,\,i,j\in\{1,...,d\}$, on a $(-1)^k\frac{\partial^k}{\partial z_j^k}\Phi^i_p(x_1,...,x_d)\geq 0$. En particulier la fonction $t\mapsto \Phi^i_p({\bf x}+te_j)$ est positive, d\'ecroissante, et convexe pour $t\geq 0$.
\end{lem}

{\bf Preuve:}  En 2 pas:
\begin{enumerate}
\item[1)] $d=1$. Comme dans la preuve de la proposition \ref{deriv}, on \'ecrit $p(z)=c\prod_{i=1}^n(z-w_i)$, o\`u les $w_i$ sont r\'eels. Alors $\Phi_p(x)=\sum_{i=1}^n \frac{1}{x-w_i}$, et donc la d\'eriv\'ee $k$-i\`eme $\Phi^{(k)}_p$ est donn\'ee par $\Phi^{(k)}_p(x)=(-1)^k \sum_{i=1}^n \frac{k!}{(x-w_i)^{k+1}}$; la derni\`ere somme est positive pour $x> ZM(p)$.
\item[2)] $d>1$. Pour $k=0$ ou $i=j$, le r\'esultat d\'ecoule du pas (1). On peut donc supposer $k\geq 1$ et $i\neq j$. Quitte \`a renum\'eroter les variables, on peut supposer $d=2$, $i=1$ et $j=2$ et il faut montrer que, si $p(z_1,z_2)$ est r\'eel stable sans z\'ero dans le quadrant $\{y_1\geq x_1,y_2\geq x_2\}$, alors
\begin{equation*}
\begin{split}
(-1)^k\frac{\partial^k}{\partial z_2^k}\Phi^1_p(x_1,x_2)\geq 0 &\Leftrightarrow (-1)^k\frac{\partial^k}{\partial z_2^k}(\frac{\partial}{\partial z_1}\log p)(x_1,x_2)\geq0 \\
& \Leftrightarrow \frac{\partial}{\partial z_1}((-1)^k\frac{\partial^k}{\partial z_2^k}\log p)(x_1,x_2)\geq 0.
\end{split}
\end{equation*}
On va donc montrer que la fonction $t\mapsto (-1)^k(\frac{\partial^k}{\partial z_2^k}\log p)(t,x_2)$ est croissante pour $t\geq x_1$. Par continuit\'e, il suffit de le montrer pour un ensemble g\'en\'erique de valeurs de $t$ (o\`u un ensemble fini de valeurs exceptionnelles est exclu).

Pour $t$ fix\'e, le polyn\^ome $q(z)=p(t,z)$ est r\'eel stable de degr\'e $n=n(t)$ (voir le lemme \ref{special}), avec des racines r\'eelles $y_1(t),...,y_n(t)$ (\'eventuellement multiples). Pour $t$ g\'en\'erique, le degr\'e $n$ de $q$ ne d\'epend pas de $t$, et on peut supposer que les racines $y_i$ d\'ependent de $t$ au moins de mani\`ere $C^1$. Comme au pas (1), on \'ecrit $q(z)=c\prod_{i=1}^n (z-y_i(t))$ donc
\begin{equation*}
\begin{split}
(-1)^k(\frac{\partial^k}{\partial z_2^k}\log p)(t,x_2) &=(-1)^k(\log q)^{(k)}|_{z=x_2} \\
&=(-1)^k(\sum_{i=1}^n \log(z-y_i(t))^{(k)}|_{z=x_2}=-\sum_{i=1}^n\frac{(k-1)!}{(x_2-y_i(t))^k},
\end{split}
\end{equation*}
et il suffit de montrer que $t\mapsto \frac{1}{x_2-y_i(t)}$ est d\'ecroissante, pour $i=1,...,n$ et $t\geq x_1$. Comme $p$ est sans z\'ero dans le quadrant $\{y_1\geq x_1,y_2\geq x_2\}$, on a $x_2>ZM(q)$, c-\`a-d. $x_2> y_i(t)$ pour $i=1,...,n$. Donc il suffit de montrer que $t\mapsto y_i(t)$ est d\'ecroissante pour $i=1,...,n$. Si ce n'\'etait pas le cas, on trouverait un indice~$i$ et une valeur $\tau_1$ de $t$ telle que $y_i'(\tau_1)>0$. Soit $\tau_2=y_i(\tau_1)$. Consid\'erons le d\'eveloppement de Taylor de $p$ autour de $(\tau_1,\tau_2)$:
$$p(z_1,z_2)=p(\tau_1,\tau_2)+(z_1-\tau_1)(\partial_1p)(\tau_1,\tau_2)+(z_2-\tau_2)(\partial_2p)(\tau_1,\tau_2)+ O(|z_1-\tau_1|^2+|z_2-\tau_2|^2).$$
Mais $p(\tau_1,\tau_2)=0$ et $(\partial_1p)(t,z)=- c\sum_{k=1}^n y_k'(t)\prod_{j\neq k}(z-y_j(t))$, donc $(\partial_1p)(\tau_1,\tau_2)=-cy_i'(\tau_1)\prod_{j\neq i}(\tau_2-y_j(\tau_1))$; tandis que 
$(\partial_2p)(t,z)=c\sum_{k=1}^n\prod_{j\neq k}(z-y_j(t))$, donc $(\partial_2p)(\tau_1,\tau_2)=c\prod_{j\neq i}(\tau_2-y_j(\tau_1))$. Ainsi:
$$p(z_1,z_2)=c\prod_{j\neq i}(\tau_2-y_j(\tau_1))[(z_2-\tau_2)-y_i'(\tau_1)(z_1-\tau_1)] + O(|z_1-\tau_1|^2+|z_2-\tau_2|^2).$$
Comme $y_i'(\tau_1)>0$, la droite complexe d'\'equation $(z_2-\tau_2)-y_i'(\tau_1)(z_1-\tau_1)=0$ rencontre l'ouvert $\mathbb{H}^2$ de $\CCC^2$. Par le th\'eor\`eme des fonctions implicites, $p$ a des z\'eros voisins de $(\tau_1,\tau_2)$ dans $\mathbb{H}^2$, ce qui contredit le caract\`ere r\'eel stable de $p$.
\hfill$\square$
\end{enumerate}

\begin{lem}[\cite{MSS-KS}, lemma 5.9; \cite{Tao}, lemma 20]
\label{lem19+20}
Soit ${\bf x}=(x_1,...,x_d)\in\RRR^d$. Soit $p(z_1,...,z_d)$ un polyn\^ome r\'eel stable, sans z\'ero dans l'orthant $\{{\bf y\geq x}\}$. On suppose que, pour un certain $j\in\{1,...,d\}$, il existe $\delta>0$ tel que $\Phi^j_p(x_1,...,x_d)+\frac{1}{\delta}\leq 1$. Alors $(1-\partial_j)p$ n'a pas de z\'ero dans l'orthant $\{{\bf y\geq x+\delta e_j}\}$, et de plus pour tout $i=1,...,d$:
\begin{equation}\label{estimate}
\Phi^i_{(1-\partial_j)p}({\bf x+\delta e_j})\leq \Phi^i_p({\bf x}).
\end{equation}
\end{lem}

{\bf Preuve:} Si ${\bf y}$ est dans l'orthant $\{{\bf y\geq x+\delta e_j}\}$, alors par le lemme \ref{lem16+17} on a $\Phi^j_p({\bf y})\leq \Phi^j_p({\bf x})\leq 1-\frac{1}{\delta}<1$, donc $\partial_jp({\bf y})<p({\bf y})$ et ainsi $(1-\partial_j)p({\bf y})>0$.

Pour montrer l'in\'egalit\'e (\ref{estimate}), partons de $(1-\partial_j)p= p(1-\Phi^j_p)$ et prenons la d\'eriv\'ee logarithmique par rapport \`a la $i$-i\`eme variable: $\Phi^i_{(1-\partial_j)p}=\Phi^i_p-\frac{\partial_i\Phi^j_p}{1-\Phi^j_p}$. Donc (\ref{estimate}) est \'equivalente \`a: 
\begin{equation}\label{technical}
-\frac{\partial_i\Phi^j_p({\bf x+\delta e_j})}{1-\Phi^j_p({\bf x+\delta e_j})}\leq \Phi^i_p({\bf x})-\Phi^i_p({\bf x+\delta e_j}).
\end{equation}

Pour montrer (\ref{technical}), on remarque d'abord que, par le lemme \ref{lem16+17}, on a $\Phi^j_p({\bf x+\delta e_j})\leq \Phi^j_p({\bf x})$ et, par hypoth\`ese, $\Phi^j_p(x_1,...,x_d)\leq 1-\frac{1}{\delta}$, donc $\Phi^j_p({\bf x+\delta e_j})\leq 1-\frac{1}{\delta}$ d'o\`u:
$$ \delta \geq \frac{1}{1-\Phi^j_p({\bf x+\delta e_j})}.$$
Comme $p$ est sans z\'ero dans l'orthant $\{{\bf y\geq x+\delta e_j}\}$, le lemme \ref{lem16+17} (appliqu\'e en ${\bf x+\delta e_j}$) donne $-\partial_i\Phi^j_p({\bf x +\delta e_j})\geq 0$, donc 
$$-\delta.\partial_i\Phi^j_p({\bf x +\delta e_j})\geq -\frac{\partial_i\Phi^j_p({\bf x+\delta e_j})}{1-\Phi^j_p({\bf x+\delta e_j})}.$$
Pour \'etablir (\ref{technical}), il suffit donc de montrer $-\delta.\partial_i\Phi^j_p({\bf x +\delta e_j})\leq \Phi^i_p({\bf x})-\Phi^i_p({\bf x+\delta e_j})$.
Mais 
$$\partial_i\Phi^j_p({\bf x +\delta e_j})=(\partial_i\partial_j\log p)({\bf x+\delta e_j})=(\partial_j\partial_i\log p)({\bf x +\delta e_j})=\partial_j\Phi^i_p({\bf x +\delta e_j}),$$
donc on veut d\'emontrer
$$-\delta.\partial_j\Phi^i_p({\bf x +\delta e_j})\leq \Phi^i_p({\bf x})-\Phi^i_p({\bf x+\delta e_j})\Longleftrightarrow \Phi^i_p({\bf x+\delta e_j})\leq \Phi^i_p({\bf x})+\delta.\partial_j\Phi^i_p({\bf x +\delta e_j})$$
qui exprime exactement la convexit\'e de la fonction $t\mapsto \Phi^i_p({\bf x}+te_j)$, \'etablie au lemme \ref{lem16+17}.
\hfill$\square$

\begin{prop}[\cite{Tao}, Cor. 21] 
\label{Cor21}
Soit ${\bf x}=(x_1,...,x_d)\in\RRR^d$. Soit $p(z_1,...,z_d)$ un polyn\^ome r\'eel stable, sans z\'ero dans l'orthant $\{{\bf y\geq x}\}$. S'il existe $\delta>0$ tel que $\Phi^j_p(x_1,...,x_d)+\frac{1}{\delta}\leq 1$ (pour tout $j=1,...,d$), alors le polyn\^ome $(\prod_{i=1}^d (1-\partial_i))p$ est sans z\'ero dans l'orthant $\{{\bf y\geq x+d}\}$, o\`u ${\bf d}=(\delta,...,\delta)$.
\end{prop}

{\bf Preuve:} Pour $0\leq k\leq d$, posons ${\bf x^{(k)}}=(x_1+\delta,...,x_k+\delta,x_{k+1},..., x_d)$, et \mbox{$q_k=(\prod_{i=1}^k (1-\partial_i))p$.} Par r\'ecurrence sur $k$ \`a partir du lemme \ref{lem19+20}, on voit que $q_k$ n'a pas de z\'ero dans l'orthant $\{{\bf y\geq x^{(k)}}\}$ et $\Phi^j_{q_k}({\bf x^{(k)}})+\frac{1}{\delta}\leq 1$ pour $j=1,...,d$. On obtient le r\'esultat d\'esir\'e pour $k=d$.
\hfill$\square$

\section{Polyn\^omes caract\'eristiques mixtes}\label{carmixte}

\subsection{D\'efinitions}

\begin{Def}\label{polmixte} Pour des matrices $A_1,...,A_d\in M_m(\CCC)$, le polyn\^ome caract\'e\-ristique mixte de la famille $\{A_1,...,A_d\}$ est:
$$\mu[A_1,...,A_d](z)= \Bigl (\prod_{i=1}^d(1-\partial_i)\Bigr ) \det \Bigl (z.1_m+\sum_{i=1}^d z_iA_i\Bigr )|_{z_1=...=z_d=0}.$$
\end{Def}

\begin{exem} Si $A_1,...,A_d$ sont des matrices semi-d\'efinies positives, il r\'esulte du lemme \ref{mainex} et des propositions \ref{special} et \ref{deriv} que $\mu[A_1,...,A_d](z)$ est un polyn\^ome r\'eel stable, c-\`a-d. a toutes ses racines r\'eelles.
\end{exem}

La principale motivation pour la d\'efinition \ref{polmixte}, vient du lemme suivant: 

\begin{lem}[\cite{Tao}, Prop. 3]
\label{Prop3}
Soient $A_1,...,A_d$ des matrices de rang $1$ dans $M_m(\CCC)$, et $ A=\sum_{i=1}^d A_i$. Alors 
$$p_A(z)=\mu[A_1,...,A_d](z).$$
\end{lem}

{\bf Preuve:} On commence par observer que, pour toute matrice $B\in M_m(\CCC)$, le \mbox{polyn\^ome} $(z_1,...,z_d)\mapsto \det(B+\sum_{i=1}^d z_iA_i)$ est {\it affine-multilin\'eaire}, c-\`a-d. de la forme 
$$\sum_{1\leq i_1<i_2<...<i_j\leq d} c_{i_1,...,i_d} z_{i_1}...z_{i_j}$$
(ou encore: dans chaque mon\^ome, chaque variable appara\^it avec un degr\'e $\leq 1$). Le cas $d=1$ se traite en travaillant dans une base dont le premier vecteur est dans l'image 
de~$A_1$, et en d\'eveloppant le d\'eterminant par rapport \`a la premi\`ere ligne. Le cas g\'en\'eral se ram\`ene \`a $d=1$ en {\oge gelant\feg} $d-1$ variables.

Comme un polyn\^ome affine-multilin\'eaire est \'egal \`a son d\'eveloppement de Taylor \`a l'ordre $(1,...,1)$, on a: 
$$\det\Bigl (B+\sum_{i=1}^d t_iA_i\Bigr )=\Bigl (\prod_{i=1}^d(1+t_i\partial_i)\Bigr )\det \Bigl (B+\sum_{i=1}^d z_iA_i\Bigr )|_{z_1=...z_d=0}.$$
Le r\'esultat s'obtient alors avec $B=z.1_m$ et $t_1=...=t_d=-1$.
\hfill$\square$

\subsection{Preuve du th\'eor\`eme \ref{MSS-posproba}}

\begin{prop}[\cite{Tao}, Cor. 4 et 15]
\label{Cor4+15}
Soient $A_1,...,A_d$ des variables al\'eatoires ind\'ependantes \`a valeurs dans les matrices semi-d\'efinies positives de rang $1$ de $M_m(\CCC)$, et $A=\sum_{i=1}^d A_i$. Alors:
\begin{enumerate}
\item[1)] $\EE p_A(z)=\mu[\EE A_1,...,\EE A_d](z)$.
\item[2)] Supposons de plus que les $A_i$ prennent un nombre fini de valeurs. Alors, pour tout $j\in\{1,...,d\}$ et toute r\'ealisation de $A_1,...,A_{j-1}$, on a
$$ZM(\mu[A_1,...,A_{j-1},A_j,\EE A_{j+1},...,\EE A_d])\leq ZM(\mu[A_1,...,A_{j-1},\EE A_j,\EE A_{j+1},...,\EE A_d])$$
pour au moins une r\'ealisation de $A_j$.
\end{enumerate}
\end{prop}

{\bf Preuve:} \begin{enumerate}
\item[1)] Par le lemme \ref{Prop3}, on a $p_A(z)=\mu[A_1,...,A_d](z)$ et la preuve du m\^eme lemme montre que $\mu[A_1,...,A_d](z)$ est une somme de termes multilin\'eaires en les coefficients de $A_1,...,A_d$.  En utilisant l'ind\'ependance des $A_i$, on voit alors que $\EE\mu[A_1,...,A_d](z)=\mu[\EE A_1,...,\EE A_d](z)$.
\item[2)] Comme $\EE A_j$ est une combinaison convexe des valeurs de $A_j$, vu le caract\`ere affine-multilin\'eaire de $\mu$, le polyn\^ome $\mu[A_1,...,A_{j-1},\EE A_j,\EE A_{j+1},...,\EE A_d](z)$ est une combinaison convexe des polyn\^omes $\mu[A_1,...,A_{j-1},A_j,\EE A_{j+1},...,\EE A_d](z)$, lesquels sont des polyn\^omes r\'eels stables. Par le lemme \ref{lem14}, $$ZM(\mu[A_1,...,A_{j-1},\EE A_j,\EE A_{j+1},...,\EE A_d])$$ est dans l'enveloppe convexe des $ZM(\mu[A_1,...,A_{j-1},A_j,\EE A_{j+1},...,\EE A_d])$. L'\'ev\'enement d\'ecrit dans l'\'enonc\'e se produit donc pour au moins une r\'ealisation de $A_j$. 
\hfill$\square$
\end{enumerate}

{\bf Preuve du th\'eor\`eme \ref{MSS-posproba}:} Par le point (1) de la proposition \ref{Cor4+15}:
$$ZM(\EE p_A)=ZM(\mu[\EE A_1,...,\EE A_d]).$$
Puis on d\'eroule les $d$ in\'egalit\'es du point (2) de la proposition \ref{Cor4+15}; pour au moins une r\'ealisation de $A_1$:
$$ZM(\mu[\EE A_1,...,\EE A_d])\geq ZM(\mu[A_1,\EE A_2,...,\EE A_d])\,;$$
pour au moins une r\'ealisation de $A_2$:
$$ZM(\mu[A_1,\EE A_2,...,\EE A_d])\geq ZM(\mu[A_1,A_2, \EE A_3,...,\EE A_d]), $$
etc. Au final, pour au moins une r\'ealisation de $A_1,A_2,...,A_d$:
$$ZM(\EE p_A)\geq ZM(\mu[A_1,...,A_d]).$$
Enfin, par le lemme \ref{Prop3}, on a: $ZM(\mu[A_1,...,A_d])=ZM(p_A)=\|A\|$.
\hfill$\square$

\subsection{Preuve du th\'eor\`eme \ref{MSS-expect}}

\begin{prop}[\cite{MSS-KS}, Thm. 5.1; \cite{Tao}, Thm. 18]
\label{Thm18}
Soient $A_1,...,A_d\in M_m(\CCC)$ des matrices semi-d\'efinies positives, avec $\sum_{i=1}^d A_i=1_m$ et $tr(A_i)\leq \varepsilon$ pour $i=1,...,d$. Soit $p(z_1,...,z_d)= \det(\sum_{i=1}^d z_iA_i)$. Alors le le polyn\^ome $(\prod_{i=1}^d(1-\partial_i))p$ n'a pas de z\'ero dans l'orthant $\{{\bf y\geq e}\}$, o\`u ${\bf e}=((1+\sqrt{\varepsilon})^2,..., (1+\sqrt{\varepsilon})^2)$.
\end{prop}

{\bf Preuve:} On veut appliquer la proposition \ref{Cor21} au polyn\^ome $p$. Celui-ci est r\'eel stable, comme sp\'ecialisation du polyn\^ome $(z,z_1,...,z_d)\mapsto\det(z.1_m+\sum_{i=1}^d z_iA_i)$ (voir le lemme \ref{mainex} et la proposition \ref{special}). Soient $t>0$ un param\`etre positif et ${\bf t}=(t,...,t)$. Montrons que $p$ est sans z\'ero dans l'orthant $\{{\bf x\geq t}\}$. En effet, pour $(x_1,...,x_d)$ dans cet orthant, on a $\sum_{i=1}^d x_iA_i\geq \sum_{i=1}^d tA_i=t.1_m$, donc $\sum_{i=1}^d x_iA_i$ est inversible. Il reste \`a v\'erifier la condition sur les $\Phi^j_p$. On va pour cela utiliser la {\it formule de Jacobi:} si $t\mapsto A(t)$ est une fonction d\'erivable \`a valeurs dans $M_m(\CCC)$, on a $(\det\,A(t))'= tr(adj(A(t)).A'(t))$, o\`u $adj$ d\'esigne la matrice des co-facteurs; en particulier, si $A(t)$ est inversible, la d\'eriv\'ee logarithmique de $\det\,A(t)$ est donn\'ee par: $\frac{(\det\,A(t))'}{\det\,A(t)}=tr(A(t)^{-1}.A'(t))$.

Ici, avec $(x_1,...,x_d)$ dans l'orthant strictement positif: 
$$\Phi^j_p(x_1,...,x_d)= tr\Bigl (\Bigl (\sum_{i=1}^d x_iA_i\Bigr  )^{-1}.A_j\Bigr ).$$
En particulier $\Phi^j_p(t,...,t)=\frac{tr(A_j)}{t}\leq\frac{\varepsilon}{t}$. Prenons alors $t=\varepsilon +\sqrt{\varepsilon}$ et $\delta=1+\sqrt{\varepsilon}$, de sorte que $\frac{\varepsilon}{t}+\frac{1}{\delta}=1$, et donc $\Phi^j_p(t,...,t)+\frac{1}{\delta}\leq 1$. En remarquant que $t+\delta=(1+\sqrt{\varepsilon})^2$, la Proposition \ref{Cor21} s'applique et donne le r\'esultat d\'esir\'e.
\hfill$\square$.

\medskip
{\bf Preuve du th\'eor\`eme \ref{MSS-expect}:} Par le point (1) de la proposition \ref{Cor4+15}, on a:
$$\EE p_A(z)=\Bigl (\prod_{i=1}^d (1-\partial_i)\Bigr )\det\Bigl (z.1_m+\sum_{i=1}^d z_i\EE A_i\Bigr )|_{z_1=...=z_d=0}.$$
Mais $z.1_m+\sum_{i=1}^d z_i\EE A_i=\sum_{i=1}^d (z+z_i)\EE A_i$, donc par une translation de $z$ sur les variables $z_i$ on a:
$$\EE p_A(z)=\Bigl (\prod_{i=1}^d (1-\partial_i)\Bigr )\det\Bigl (\sum_{i=1}^d z_i\EE A_i\Bigr )|_{z_1=...=z_d=z}.$$
En notant que $tr(\EE A_i)=\EE(tr(A_i))=\EE\|A_i\|\leq \varepsilon$ (car $A_i$ est semi-d\'efinie positive de rang 1), on peut appliquer la proposition \ref{Thm18} aux matrices $\EE A_i$: le polyn\^ome \mbox{$(\prod_{i=1}^d (1-\partial_i))\det(\sum_{i=1}^d z_i\EE A_i)$} n'a pas de z\'ero dans l'orthant $\{{\bf y\geq e}\}$, o\`u 
\mbox{${\bf e}=((1+\sqrt{\varepsilon})^2,..., (1+\sqrt{\varepsilon})^2)$,} donc $ZM(\EE p_A)\leq (1+\sqrt{\varepsilon})^2$.
\hfill$\square$

\section{Solution du probl\`eme de Kadison-Singer}\label{finKS}

Dans \cite{Wea}, N. Weaver introduit pour chaque entier $r\geq 2$ une conjecture $KS_r$, et montre que chacune est \'equivalente au probl\`eme de Kadison-Singer:

\medskip
$(KS_r)$: Soit $r\geq 2$. Il existe des constantes universelles $N\geq 2$ et $\varepsilon>0$ telles que~: si $A_1,...,A_d\in M_m(\CCC)$ sont des matrices semi-d\'efinies positives de rang 1, avec $\sum_{i=1}^d A_i=N.1_m$ et $\|A_i\|\leq 1$ pour $i=1,...,d$, alors il existe une partition $\{S_1,...,S_r\}$ de $\{1,...,d\}$ telle que $\|\sum_{i\in S_j} A_i\|\leq N-\varepsilon$.

\medskip
Au corollaire 1.3 de \cite{MSS-KS}, Marcus, Spielman et Srivastava d\'emontrent la conjecture $KS_2$ avec $N=18$ et $\varepsilon=2$. Nous adoptons la forme l\'eg\`erement diff\'erente propos\'ee par T. Tao (\cite{Tao}, Theorem 22).

\begin{thm}\label{Weaver} On fixe des entiers $d,m,r\geq 2$ et une constante $C>0$. Soient $A_1,...,A_d\in M_m(\CCC)$ des matrices semi-d\'efinies positives de rang $1$, avec $\|A_i\|\leq C$ pour $i=1,...,d$ et $\sum_{i=1}^d A_i=1_m$. Il existe une partition $\{S_1,...,S_r\}$ de $\{1,...,d\}$ telle que $\|\sum_{i\in S_j} A_i\|\leq (\sqrt{\frac{1}{r}}+\sqrt{C})^2$ pour $j=1,...,r$.
\end{thm}

{\bf Preuve:} Pour $i=1,...,d$, on note $E_i$ la variable al\'eatoire uniforme sur $\{1,...,r\}$, qui consiste \`a tirer au hasard, avec probabilit\'e $\frac{1}{r}$, un projecteur $P_j$ sur un des $r$ vecteurs de base $e_j$ de $\ell^2_r$. Pour $i=1,...,d$, on pose alors $\tilde{A}_i=r(A_i\otimes E_i)$, et $\tilde{A}=\sum_{i=1}^d \tilde{A}_i$. On a alors:
\begin{equation*}
\begin{split}
\EE\tilde{A} = r\sum_{i=1}^d \EE(A_i\otimes E_i) &=\sum_{i=1}^d \sum_{j=1}^r (A_i\otimes P_j) \\
&=\sum_{j=1}^r(\sum_{i=1}^d A_i)\otimes P_j =\sum_{j=1}^r (1_m\otimes P_j)=1_m\otimes(\sum_{j=1}^r P_j)=1_m\otimes 1_r.
\end{split}
\end{equation*}
De plus: $\|\tilde{A}_i\|=r\|A_i\|\leq rC$.

On applique le th\'eor\`eme \ref{MSS-posproba}, puis le th\'eor\`eme \ref{MSS-expect} (avec $\varepsilon=rC$): pour au moins une r\'ealisation des $\tilde{A}_i$:
$$\|\tilde{A}\|\leq ZM(\EE p_{\tilde{A}})\leq (1+\sqrt{rC})^2.$$
Cette r\'ealisation des $\tilde{A}_i$ permet de d\'efinir la partition $\{S_1,...,S_r\}$ de $\{1,...,d\}$:
$$S_j=\{i\in\{1,...,d\}: E_i=P_j\}\;\;(j=1,...,r).$$
On a alors: $$r\|\sum_{i\in S_j}A_i\|=r\|\sum_{i\in S_j}(A_i\otimes P_j)\|=\|\sum_{i\in S_j}\tilde{A}_i\|\leq \|\sum_{i=1}^d \tilde{A}_i\|=\|\tilde{A}\|
\leq (1+\sqrt{rC})^2.$$
On a alors le r\'esultat en divisant par $r$.
\hfill$\square$

\begin{lem}[\cite{Tao}, Cor. 23]
\label{Cor23}
 Soit $P\in M_d(\CCC)$ un projecteur. Pour tout $r\in\NN$, il existe des projecteurs diagonaux $Q_1,...,Q_r\in M_d(\CCC)$ avec $\sum_{i=1}^r Q_i=1_d$ et 
$$\|Q_iPQ_i\|\leq (\sqrt{\frac{1}{r}}+\sqrt{\|diag(P)\|_\infty})^2.$$
\end{lem}

{\bf Preuve:} Notons $V\subset\ell^2_d$ l'image de $P$, et $m=\dim V$ son rang. On d\'efinit une famille $A_1,...,A_d$ d'op\'erateurs de rang 1 sur $V$ par: $A_i(v)=\langle v|P(e_i)\rangle P(e_i)$. Alors $\|A_i\|=\|P(e_i)\|^2_2=\langle P(e_i)|e_i\rangle\leq\|diag(P)\|_\infty$. De plus $\sum_{i=1}^d A_i=1_m$ car 
\begin{equation*}
\begin{split}
\Bigl\langle \Bigl (\sum_{i=1}^d A_i\Bigr )v|v\Bigr \rangle &= \sum_{i=1}^d \langle A_i(v)|v\rangle = \sum_{i=1}^d |\langle v|P(e_i)\rangle|^2\\
&= \sum_{i=1}^d |\langle P(v)|e_i\rangle|^2=\sum_{i=1}^d |\langle v|e_i\rangle|^2=\|v\|_2^2
\end{split}
\end{equation*}
puisque $v\in V$. On applique le th\'eor\`eme \ref{Weaver} avec $C=\|diag(P)\|_\infty$: il existe une partition $\{S_1,...,S_r\}$ de $\{1,...,d\}$ avec 
$$\|\sum_{i\in S_j} A_i\|\leq (\sqrt{\frac{1}{r}}+\sqrt{\|diag(P)\|_\infty})^2,$$
pour $j=1,...,r$. Notons alors $Q_j$ le projecteur diagonal de $M_d(\CCC)$ correspondant \`a $S_j$. On a, pour $j=1,...,r$:
$$\|Q_jPQ_j\|=\|(Q_jP)(Q_jP)^*\|=\|Q_jP\|^2=\|Q_j|_V\|^2.$$
Mais, pour $v\in V$ (donc $v=P(v)$):
\begin{equation*}
\begin{split}
\|Q_j(v)\|_2^2&=\sum_{i\in S_j}|\langle v|e_i\rangle|^2=\sum_{i\in S_j} |\langle P(v)|e_i\rangle|^2=\sum_{i\in S_j} |\langle v|P(e_i)\rangle|^2\\
&=\sum_{i\in S_j} \langle A_i(v)|v\rangle\leq\|\sum_{i\in S_j}A_i\|\|v\|_2^2\leq(\sqrt{\frac{1}{r}}+\sqrt{\|diag(P)\|_\infty})^2\|v\|_2^2,
\end{split}
\end{equation*}
ce qui termine la preuve.
\hfill$\square$

\medskip
L'\'equivalence de l'\'enonc\'e suivant avec le probl\`eme de Kadison-Singer appara\^it dans l'article original de Kadison-Singer \cite[lemma 5]{KS}. 

\begin{thm}\label{Anderson} Pour tout $\varepsilon>0$, il existe $r\in\NN$ tel que, pour tout $T\in{\mathcal B}(\ell^2)$ avec $diag(T)=0$, il existe une famille $Q_1,...,Q_r$ de projecteurs diagonaux, avec $\sum_{i=1}^r Q_i=1$ et $\|Q_iTQ_i\|\leq\varepsilon \|T\|$ pour $i=1,...,r$.
\end{thm}

Ce r\'esultat affirme donc que, pour tout op\'erateur $T$ de diagonale nulle, on peut trouver une partition finie $\{A_1,...,A_r\}$ de $\NN$, telle que la norme des op\'erateurs compress\'es $(T_{mn})_{m,n\in A_i}$ (pour $i=1,...,r$) soit arbitrairement petite. L'\'etude des op\'erateurs de diagonale nulle est naturelle, en vue du probl\`eme de Kadison-Singer.

\begin{exem} Consid\'erons l'op\'erateur $S$ de d\'ecalage unilat\'eral sur $\ell^2$, d\'efini par $S{\bf e}_n={\bf e}_{n+1}$ pour tout $n\in\NN$. Cet op\'erateur v\'erifie la conjecture de pavage avec $\varepsilon=0$: il suffit de prendre pour $P_1$ l'ensemble $2\NN+1$ des nombres impairs, et pour $P_2$ l'ensemble $2\NN$ des nombres pairs; on a clairement $Q_1SQ_1=Q_2SQ_2=0$.
\end{exem}

\medskip
{\bf Preuve du th\'eor\`eme \ref{Anderson}:} La preuve se fait en trois pas: on traite d'abord le cas des matrices auto-adjointes de taille finie; puis on passe aux matrices de taille finie quelconques; enfin un argument de compacit\'e permet de passer de la dimension finie \`a la dimension infinie. L'\'enonc\'e pour des matrices de taille finie est connu dans la litt\'erature sous le nom de {\it conjecture de pavage de J. Anderson}; pour l'\'equivalence avec le probl\`eme de Kadison-Singer, voir \cite[Theorem 3.6]{AndTAMS}.

\smallskip
\underline{Premier pas} (d'apr\`es une id\'ee de P. Casazza): Soit $T=T^*\in M_d(\CCC)$ avec 
\mbox{$diag(T)=0$}; on peut clairement supposer $\|T\|\leq 1$. On forme alors, dans $M_{2d}(\CCC)$, la matrice:
$$P= \left(\begin{array}{cc}\frac{1+T}{2} & \frac{1}{2}\sqrt{1-T^2} \\ \frac{1}{2}\sqrt{1-T^2} & \frac{1-T}{2}\end{array}\right);$$
c'est un projecteur qui de plus v\'erifie $diag(P_1)=\bigl (\frac{1}{2},...,\frac{1}{2}\bigr )$. On prend $r$ assez grand pour avoir $2\Bigl (\sqrt{\frac{1}{r}}+\sqrt{\frac{1}{2}}\Bigr )^2-1 \leq \varepsilon$. Par le lemme \ref{Cor23}, on trouve des projecteurs diagonaux $Q''_1,...,Q''_r\in M_{2d}(\CCC)$, avec $\sum_{i=1}^r Q''_i=1_{2d}$ et $\|Q''_iPQ''_i\|\leq \Bigl (\sqrt{\frac{1}{r}}+\sqrt{\frac{1}{2}}\Bigr )^2$ pour $i=1,..,r$. Dans $M_{2d}(\CCC)$, notons $E$ (resp. $E'$) le projecteur sur le sous-espace engendr\'e par les $d$~premiers (resp. derniers) vecteurs de base, de sorte que 
\mbox{$EPE=\frac{1+T}{2}$} (resp. $E'PE'=\frac{1-T}{2}$). Pour $i=1,...,r$, posons $Q_i=EQ''_i$ (resp. $Q'_i=E'Q''_i$), donc 
\mbox{$\sum_{i=1}^rQ_i=\sum_{j=1}^rQ'_j=1_d$} et, pour $i=1,...,d$:
$$\|Q_i(1+T)Q_i\|\leq 2\Bigl (\sqrt{\frac{1}{r}}+\sqrt{\frac{1}{2}}\Bigr )^2\;(\mbox{resp.}\;\|Q'_i(1-T)Q'_i\|\leq 2\Bigl (\sqrt{\frac{1}{r}}+\sqrt{\frac{1}{2}}\Bigr )^2).$$
Donc $$-Q_i\leq Q_iTQ_i\leq \Bigl [2\Bigl (\sqrt{\frac{1}{r}}+\sqrt{\frac{1}{2}}\Bigr )^2-1\Bigr ]Q_i\leq\varepsilon Q_i$$
 (resp. $$-\varepsilon Q'_i\leq \Bigl [1-2\Bigl (\sqrt{\frac{1}{r}}+\sqrt{\frac{1}{2}}\Bigr )^2\Bigr ]Q'_i\leq Q'_iTQ'_i\leq Q'_i).$$
  Posons $Q_{i,j}=Q_iQ'_j$ ($1\leq i,j\leq r$). On a alors $\sum_{i,j=1}^r Q_{i,j}=1_d$ et:
  $$-\varepsilon Q_{i,j}\leq Q_{i,j}TQ_{i,j}\leq \varepsilon Q_{i,j},$$
  donc $\|Q_{i,j}TQ_{i,j}\|\leq\varepsilon$.

\medskip
\underline{Deuxi\`eme pas}: Soit $T\in M_d(\CCC)$ avec $diag(T)=0$. On \'ecrit $T=\frac{T+T^*}{2} + i(\frac{T-T^*}{2i})$, et on applique le premier cas aux matrices auto-adjointes $A=\frac{T+T^*}{2}$ et $B=\frac{T-T^*}{2i}$ (en observant que $\|A\|\leq\|T\|$ et $\|B\|\leq\|T\|$). On trouve donc des projecteurs diagonaux $Q'_1,...,Q'_r,Q''_1,...,Q''_r\in M_d(\CCC)$ avec $\sum_{i=1}^r Q'_i=\sum_{j=1}^r Q''_j=1_d$ et $\|Q'_iAQ'_i\|\leq\frac{\varepsilon}{2}\|T\|$ pour $i=1,...,r$ et $\|Q''_jBQ''_j\|\leq\frac{\varepsilon}{2}\|T\|$ pour $j=1,...,r$. On pose alors $Q_{ij}=Q'_iQ''_j$ et on a $\sum_{i,j} Q_{ij}=1_d$ avec $\|Q_{ij}TQ_{ij}\|\leq\varepsilon\|T\|$ pour $1\leq i,j\leq r$.

\medskip
\underline{Troisi\`eme pas:} Soit $T\in{\mathcal B}(\ell^2)$, avec $diag(T)=0$. Notons $E_d:\ell^2\rightarrow \ell^2_d$ le projecteur (diagonal), et $T_d=E_dTE_d$. Par le deuxi\`eme pas, on trouve des projecteurs diagonaux $Q^{(d)}_1,...,Q^{(d)}_r$ avec $\sum_{i=1}^r Q^{(d)}_i=E_d$ et $\|Q_i^{(d)}T_dQ_i^{(d)}\|\leq\varepsilon\|T\|$ pour $i=1,...,r$. On identifie l'ensemble des projecteurs diagonaux \`a l'espace compact $\{0,1\}^\NN$. Un argument diagonal, appliqu\'e aux $r$ suites $(Q_i^{(d)})_{d\geq 1}$ de $\{0,1\}^\NN$, fournit une suite strictement croissante d'entiers $(d_k)_{k>0}$ et des projecteurs diagonaux $Q_1,..., Q_r$ avec $\lim_{k\rightarrow\infty} Q^{(d_k)}_i=Q_i$ pour $i=1,..,r$. On a d'une part $\sum_{i=1}^r Q_i=\lim_{k\rightarrow\infty}(\sum_{i=1}^r Q_i^{(d_k)})=\lim_{k\rightarrow\infty} 1_{d_k}=1$; d'autre part, pour $\xi,\eta\ell^2$ \`a support fini, on a $\xi,\eta\in\ell^2_{d_k}$ pour $k\gg 0$, donc:
\begin{equation*}
\begin{split}
|\langle Q_iTQ_i\xi|\eta\rangle| &= |\langle TQ_i\xi|Q_i\eta\rangle|=|\langle T_{d_k}Q_i^{(d_k)}\xi|Q_i^{(d_k)}\eta\rangle|\;\mbox{(pour}\;k\gg 0)\\
&=|\langle Q_i^{(d_k)}T_{d_k}Q_i^{(d_k)}\xi|\eta\rangle|\leq \|Q_i^{(d_k)}T_{d_k}Q_i^{(d_k)}\|\|\xi\|_2\|\eta\|_2\leq\varepsilon\|T\|\|\xi\|_2\|\eta\|_2,
\end{split}
\end{equation*}
d'o\`u $\|Q_iTQ_i\|\leq\varepsilon\|T\|$.
\hfill$\square$

\medskip
{\bf Preuve du th\'eor\`eme \ref{KS-main}:} Soit $\varphi$ un \'etat pur sur $D$, et $\psi$ un \'etat sur ${\mathcal B}(\ell^2)$ avec $\psi|_D=\varphi$. On doit montrer que $\psi(T)=\varphi(diag(T))$ pour tout $T\in{\mathcal B}(\ell^2)$. Quitte \`a remplacer $T$ par $T-diag(T)$ et \`a normer, on peut supposer que $diag(T)=0$ et $\|T\|\leq 1$, et on doit montrer que $\psi(T)=0$. 

Soit $\varepsilon>0$. On va montrer que $|\psi(T)|\leq \varepsilon $. Par le th\'eor\`eme \ref{Anderson}, on trouve des projecteurs diagonaux $Q_1,...,Q_r$ avec $\sum_{i=1}^r Q_i=1$ et $\|Q_i TQ_i\|\leq\varepsilon$. Comme $\varphi$ est multiplicatif sur $D$, on a $\psi(Q_i)=\varphi(Q_i)\in\{0,1\}$, et donc il existe un unique indice $i_0\in\{1,...,r\}$ tel que $\psi(Q_{i_0})=1$. On a alors:
$$\psi(T)=\psi((\sum_i Q_i)T(\sum_j Q_j)) = \sum_{i,j} \psi(Q_iTQ_j).$$

Montrons que, dans cette somme double, tous les termes sont nuls sauf celui correspondant \`a $i=j=i_0$. En effet, comme le produit scalaire $(A,B)\mapsto \psi(AB^*)$ est semi-d\'efini positif sur ${\mathcal B}(\ell^2)$, on a par l'in\'egalit\'e de Cauchy-Schwarz: $|\psi(Q_iTQ_j)|^2 \leq \psi((Q_iT)(Q_iT)^*)\psi(Q_jQ_j^*)$ et donc $\psi(Q_iTQ_j)=0$ si $j\neq i_0$. Sym\'etriquement, on a $\psi(Q_iTQ_j)=0$ si $i\neq i_0$. On a donc finalement:
\mbox{$|\psi(T)|=|\psi(Q_{i_0}TQ_{i_0})|\leq\|Q_{i_0}TQ_{i_0}\|\leq\varepsilon$,} ce qui termine la preuve.
\hfill$\square$

\section{Conjectures de Bourgain-Tzafriri et de Feichtinger}\label{BTF}

Cette section doit beaucoup aux notes de cours de P. Casazza \cite{Cas}.

\begin{Def} Soit $(\varphi_n)_{n\in \NN}$ une suite dans un espace de Hilbert $\mathcal{H}$.
\begin{itemize}
\item La suite $(\varphi_n)_{n\in \NN}$ est une {\bf suite de Riesz} s'il existe des constantes $A,B>0$ telles que, pour tout ${\bf a}=(a_n)_{n\in\NN}\in\ell^2$:
$$A\|{\bf a}\|^2_2\leq \|\sum_{n\in\NN} a_n\varphi_n\|^2\leq B\|{\bf a}\|_2^2.$$
Si $A=1-\varepsilon$ et $B=1+\varepsilon$, on parle de $\varepsilon$-suite de Riesz.
\item La suite $(\varphi_n)_{n\in \NN}$ est une {\bf frame}, ou {\bf structure oblique}, s'il existe des constantes $A,B>0$ telles que, pour tout $\xi\in\mathcal{H}$:
$$A\|\xi\|^2\leq\sum_{n\in\NN}|\langle\xi|\varphi_n\rangle|^2\leq B\|\xi\|^2.$$
\end{itemize}
\end{Def}

La proposition suivante confirme une conjecture de Casazza et Vershynin (appel\'ee conjecture $R_\varepsilon$ dans \cite{Cas}):

\begin{prop}\label{Repsilon} Soit $\varepsilon >0$. Toute suite de Riesz de $\ell^2$ form\'ee de vecteurs de norme $1$, est la r\'eunion disjointe d'un nombre fini de $\varepsilon$-suites de Riesz.
\end{prop}

{\bf Preuve:} Soit $(\varphi_n)_{n\in\NN}$ une suite de Riesz de $\ell^2$, de constantes $A,B$, avec $\|\varphi_n\|_2=1$ pour tout $n$. Soit $T\in\mathcal{B}$$(\ell^2)$ l'op\'erateur born\'e d\'efini par $t({\bf e}_n)=\varphi_n$, pour tout $n\in \NN$. Posons $S=T^*T$; vu les hypoth\`eses, on a $diag(S)=(1,1,1,...)$. On peut donc appliquer le th\'eor\`eme \ref{Anderson} (conjecture de pavage) \`a l'op\'erateur $S-1$: il existe une partition finie $\{A_1,...,A_r\}$ de $\NN$, telle que $\|Q_{A_j}(S-1)Q_{A_j}\|\leq \delta\|S-1\|$ pour $j=1,...,r$, o\`u $\delta=\frac{\varepsilon}{1+\|S\|}$ et $Q_{A_j}$ est le projecteur diagonal associ\'e \`a $A_j$. Notons que $r=r(\varepsilon,\|T\|)$.

Montrons que la suite $(\varphi_n)_{n\in A_j}$ est une $\varepsilon$-suite de Riesz. En effet, pour ${\bf a}\in\ell^2$:
$$\|\sum_{n\in A_j} a_n\varphi_n\|_2^2=\|T(\sum_{n\in A_j} a_n{\bf e}_n)\|_2^2=\|TQ_{A_j}{\bf a}\|_2^2=\langle Q_{A_j}SQ_{A_j}{\bf a}|{\bf a}\rangle.$$
Mais
\begin{equation*}
\begin{split}
\langle Q_{A_j}SQ_{A_j}{\bf a}|{\bf a}\rangle &=\langle Q_{A_j}{\bf a}|Q_{A_j}{\bf a}\rangle-\langle (Q_{A_j}(1-S)Q_{A_j})Q_{A_j}{\bf a}|Q_{A_j}{\bf a}\rangle\\
&\geq (1-\delta\|1-S\|)\|Q_{A_j}{\bf a}\|^2_2\geq (1-\varepsilon)\|Q_{A_j}{\bf a}\|^2_2.
\end{split}
\end{equation*}
On montre de mani\`ere analogue que $\|\sum_{n\in A_j} a_n\varphi_n\|_2^2\leq (1+\varepsilon)\|Q_{A_j}{\bf a}\|^2_2.$
\hfill$\square$

\medskip
En 2004, H. Feichtinger \'emet la conjecture suivante, popularis\'ee par P. Casazza et ses collaborateurs (voir par exemple \cite{CCLV}):

\begin{thm}\label{feichtinger} Toute frame dans $\ell^2$ consistant en vecteurs de norme $1$, est r\'eunion disjointe d'une famille finie de suites de Riesz.
\end{thm}

{\bf Preuve:} Soit $(\varphi_n)_{n\in\NN}$ une frame dans $\ell^2$, de constantes $A,B$, avec $\|\varphi_n\|_2=1$ pour tout $n$. On consid\`ere l'op\'erateur born\'e $T:\ell^2\rightarrow\ell^2:\xi\mapsto (\langle\xi|\varphi_n\rangle)_{n\in\NN}$. Son adjoint satisfait $T^*({\bf e}_n)=\varphi_n$ pour tout $n\in\NN$. Posons $S=TT^*$: comme les $\varphi_n$ sont de norme~1, on a $diag(S)=(1,1,1,...)$; exactement comme dans la preuve de la proposition \ref{Repsilon}, on peut appliquer \`a $S-1$ la conjecture de pavage (Th\'eor\`eme \ref{Anderson}) avec $\varepsilon=\frac{1}{2(1+B)}$: on trouve une partition $\{A_1,...,A_r\}$ de $\NN$ avec $\|Q_{A_j}(S-1)Q_{A_j}\|\leq \varepsilon\|S-1\|$ et on montre comme dans la preuve de la proposition \ref{Repsilon} que, pour $j=1,...,r$:
$$\frac{1}{2}\sum_{i\in A_j}|a_i|^2\leq\|\sum_{i\in A_j} a_i\varphi_i\|^2\leq\frac{3}{2}\sum_{i\in A_j}|a_i|^2,$$
pour tout ${\bf a}=(a_i)_{i\in\NN}\in\ell^2$.
\hfill$\square$

\medskip
Dans les articles \cite{BT1,BT2}, J. Bourgain et L.Tzafriri \'etudient -- en liaison d'abord avec des questions de g\'eom\'etrie des espaces de Banach, ensuite avec le probl\`eme de Kadison-Singer -- un probl\`eme d'{\it inversibilit\'e restreinte} de matrice: si $T\in M_n(\CCC)$, et $A\subset\{1,...,n\}$ est tel que $T$ est injectif en restriction au sous-espace vectoriel de $\CCC^n$ engendr\'e par les ${\bf e}_i\;(i\in A)$, de quelle taille peut-on choisir $A$? Ils d\'emontrent:

\begin{thm}[\cite{BT1}, Theorem 1.2]
\label{thmBT}
 Il existe une constante $c>0$ telle que, pour toute matrice $T\in M_n(\CCC)$ ayant des colonnes de norme $1$, il existe une partie 
\mbox{$A\subset\{1,...,n\}$} de cardinal au moins $c\,.\,\frac{n}{\|T\|^2}$, telle que $\|\sum_{i\in A}a_iT({\bf e}_i)\|_2^2\geq c^2\sum_{i\in A}|a_i|^2$ pour tout choix de nombres complexes $(a_i)_{i\in A}$.
\hfill$\square$
\end{thm}

Les cas extr\^emes $T=1_n$ (de norme 1) et $T({\bf e}_i)={\bf e}_1$ (de norme $\sqrt{n}$), montrent que cet \'enonc\'e est essentiellement optimal quant \`a la taille de la partie $A$. Voir \cite{Naor, Sr2} pour des \'enonc\'es r\'ecents d'invertibilit\'e restreinte.

\medskip
Le th\'eor\`eme \ref{thmBT} a motiv\'e l'\'enonc\'e suivant; ce dernier semble avoir longtemps appartenu \`a la tradition orale, avant d'avoir \'et\'e consign\'e par Casazza\footnote{P. Casazza, communication personnelle, 26 d\'ecembre 2013.} sous le nom de {\it conjecture de Bourgain-Tzafriri}; voir par exemple \cite{CaTr}.

\begin{thm}\label{BT} Pour tout $B > 1$, il existe un entier $r = r(B)$ qui v\'erifie: Pour tout $n$ et toute matrice $T\in M_n(\CCC)$ avec des colonnes de norme $1$ et $\|T\|\leq B$, il existe une partition $\{A_1,...,A_r\}$ de $\{1,2,...,n\}$ telle que, pour tout $j = 1, 2,..., r $ et tout choix de nombres complexes $(a_i)_{i\in A_j}$, on a:
$$\Bigl \|\sum_{i\in A_j} a_i T({\bf e}_i)\Bigr \|_2^2\geq \frac{7}{16}\sum_{i\in A_j}|a_i|^2.$$
\end{thm}

{\bf Preuve:} En rempla\c cant $T$ par $T\oplus 1$, on peut supposer que $T$ est un op\'erateur sur~$\ell^2$, \`a colonnes de norme 1. Posons $\varphi_n=(\frac{3}{5}T({\bf e}_n),\frac{4}{5}{\bf e}_n)\in \ell^2\oplus\ell^2$. Les $\varphi_n$ forment une suite de Riesz de vecteurs de norme 1. On peut donc lui appliquer la proposition \ref{Repsilon} avec $\epsilon=(\frac{3}{5})^2$: il existe une partition $\{A_1,...,A_r\}$ de $\NN$ (o\`u $r$ est fonction de $\|T\|$) telle que, pour $j=1,...,r$:
 $$\frac{16}{25}\sum_{i\in A_j}|a_i|^2\leq \Bigl \|\sum_{i\in A_j} a_i\varphi_i\Bigr \|_2^2= \frac{16}{25}\Bigl \|\sum_{i\in A_j} a_iT({\bf e}_i)\Bigr \|^2_2+\frac{9}{25}\sum_{i\in A_j}|a_i|^2.$$
 En r\'earrangeant:
 $$\frac{7}{16}\sum_{i\in A_j}|a_i|^2\leq \Bigl \|\sum_{i\in A_j}a_iT({\bf e_i})\Bigr \|^2_2.$$
 \hfill$\square$
 
 \section{Graphes de Ramanujan de degr\'e arbitraire}\label{Raman}
 
 \subsection{Th\'eorie alg\'ebrique des graphes}
 
 Soit $X=(V,E)$ un graphe fini connexe. Notons $\sim$ la {\it relation d'adjacence} sur l'ensemble $V$ des sommets:
 $$x\sim y\Longleftrightarrow \{x,y\}\in E.$$
 Un graphe $X$ est $d$-{\it r\'egulier} si tout sommet est adjacent \`a $d$ autres sommets. La {\it matrice d'adjacence} de $X$ est:
 $$(A_X)_{xy}=\left\{\begin{array}{ccc}1 & si & x\sim y \\0 & sinon. & \end{array}\right.$$
 Le r\'esultat suivant est classique (voir par exemple \cite{God}, Th\'eor\`eme 4.2 du Chapitre 2).
 
 \begin{prop} Soit $X$ un graphe fini connexe, $d$-r\'egulier.
 \begin{enumerate}
 \item[1)] Pour toute valeur propre $\lambda$ de $A_X$, on a $|\lambda|\leq d$.
 \item[2)] $d$ est valeur propre de $A_X$, de multiplicit\'e $1$.
 \item[3)] $-d$ est valeur propre de $A_X$ si et seulement si $X$ est bi-colorable. Dans ce cas, le spectre de $A_X$ est sym\'etrique par rapport \`a $0$.
 \end{enumerate}
 \end{prop}
 
 \begin{Def} Un graphe fini connexe $d$-r\'egulier $X$ est un {\it graphe de Ramanujan} si, pour toute valeur propre $\lambda$ de $A_X$, on a $|\lambda|=d$ ou $|\lambda|\leq 2\sqrt{d-1}$.
 \end{Def}
 
 Nous renvoyons \`a \cite{Lub,DSV,Val} pour la justification de la quantit\'e $2\sqrt{d-1}$. Les familles infinies de graphes de Ramanujan $d$-r\'eguliers fournissent (pour $d\geq 3$) des familles d'expanseurs optimales du point de vue spectral. Les seules constructions connues de telles familles sont pour $d$ de la forme $1+q$, o\`u $q$ est une puissance de premier, et reposent sur de la th\'eorie des nombres difficile -- qui explique la terminologie {\oge Ramanujan\feg}. La question est pos\'ee dans \cite[Problem 10.7.3]{Lub} (voir aussi \cite{Val}, bas de la page 250) de l'existence, pour tout $d\geq 3$, de familles infinies de graphes de \mbox{Ramanujan} $d$-r\'eguliers. Dans \cite[Theorem 5.5]{MSS-Ra}, Marcus, Spielman et Srivastava d\'emontrent:
 
 \begin{thm}\label{MSS-Ram} Pour tout $d\geq 3$, il existe des familles infinies de graphes de Ramanujan $d$-r\'eguliers bi-colorables.
  \end{thm}
 
Nous verrons que la preuve est non constructive. La restriction au cas bi-colorable s'explique par le fait que les techniques de Marcus-Spielman-Srivastava leur permettent de contr\^oler la plus grande valeur propre d'une matrice, mais pas la plus petite (voir le th\'eor\`eme \ref{MSS-posproba}): l'hypoth\`ese ``bi-colorable'' permet de contr\^oler toutes les valeurs propres simultan\'ement. Pour $d$ quelconque, la question de l'existence de familles infinies de graphes de Ramanujan $d$-r\'eguliers non bi-colorables est toujours ouverte. 

\subsection{2-rel\`evements}

\begin{Def} Soit $X=(V,E)$ un graphe fini connexe. Un {\it signage} de $X$ est une fonction sym\'etrique $s:V\times V\rightarrow \{-1,0,1\}$ telle que $s(x,y)\neq 0$ si et seulement si $x\sim y$.
\end{Def}

Soit $s$ un signage d'un graphe \`a $n$ sommets. La {\it matrice d'adjacence sign\'ee} $A_X^{(s)}$ est la matrice $n$-fois-$n$ d\'efinie par:
$$(A_X^{(s)})_{xy}=s(x,y).$$
Le {\it 2-rel\`evement} de $X$ associ\'e \`a $s$ est un graphe $\tilde{X}^{(s)}$ \`a $2n$ sommets, dont l'ensemble des sommets est la r\'eunion $V_1\coprod V_2$ de 2 copies de $V$, et dont les ar\^etes sont donn\'ees par:
$$ \left\{\begin{array}{ccccccc}x_1\sim y_1 & et & x_2\sim y_2 & si  & x\sim y & et & s(x,y)=1 \\x_1\sim y_2 & et & x_2\sim y_1 & si & x\sim y & et & s(x,y)=-1\end{array}\right.$$
(o\`u, pour $i=1,2$, le sommet $x_i$ de $V_i$ correspond au sommet $x$ de $V$). On observe que:

\begin{itemize}
\item si $\tilde{X}^{(s)}$ est connexe, c'est un rev\^etement double de $X$;
\item si $X$ est $d$-r\'egulier, $\tilde{X}^{(s)}$ l'est \'egalement;
\item si $X$ est bi-colorable, $\tilde{X}^{(s)}$ l'est \'egalement.
\end{itemize}

On montre alors (voir \cite[Lemme 3.1]{BiLi}) que le spectre de la matrice d'adjacence $A_{\tilde{X}^{(s)}}$ est la r\'eunion du spectre de $A_X$ et du spectre de $A_X^{(s)}$. Bilu et Linial proposent dans \cite{BiLi} la conjecture suivante:

\begin{conj}\label{BilLin} Tout graphe $d$-r\'egulier $X$ poss\`ede un signage $s$ tel que les valeurs propres de $A_X^{(s)}$ sont inf\'erieures ou \'egales \`a $2\sqrt{d-1}$ en valeur absolue.
\end{conj}

Une cons\'equence de cette conjecture est qu'un graphe de Ramanujan $d$-r\'egulier poss\`ede un 2-rel\`evement qui est encore de Ramanujan, pour lequel on peut it\'erer la construction. Donc, cette conjecture, si elle est vraie, implique l'existence de familles infinies de graphes de Ramanujan $d$-r\'eguliers, puisqu'on peut d\'emarrer avec le graphe complet $K_{d+1}$ sur $d+1$ sommets. Marcus, Spielman et Srivastava d\'emontrent la conjecture de Bilu-Linial dans le cas bi-colorable \cite[Theorem 5.3]{MSS-Ra}:

\begin{thm}\label{MSS-sign} Tout graphe $d$-r\'egulier bi-colorable poss\`ede un signage qui satisfait la conjecture \ref{BilLin}.
\end{thm}

Puisqu'on peut d\'emarrer avec le graphe biparti complet $K_{d,d}$ sur $2d$ sommets, qui est de Ramanujan, par it\'eration de la construction des 2-rel\`evements on obtient une famille infinie de graphes de Ramanujan $d$-r\'eguliers bi-colorables, ce qui d\'emontre le th\'eor\`eme \ref{MSS-Ram}.

\subsection{Preuve du th\'eor\`eme \ref{MSS-sign}}

Soit $X$ un graphe $d$-r\'egulier bi-colorable \`a $n$ sommets et $m$ ar\^etes. Pour un signage~$s$ de $X$ et une ar\^ete $e=\{u,v\}\in E$, on d\'efinit un op\'erateur $A_e^{(s)}$ semi-d\'efini positif de rang 1 sur $\ell^2(V)$:
$$A_e^{(s)}(f)=\left\{\begin{array}{ccc}\langle f|\delta_u-\delta_v\rangle(\delta_u-\delta_v) & si & s(u,v)=-1 \\\langle f|\delta_u+\delta_v\rangle(\delta_u+\delta_v) & si & s(u,v)=1\end{array}\right.$$
($f\in\ell^2(V)$). Matriciellement, si $s(u,v)=-1$:
$$(A_e^{(s)})_{xy}=\left\{\begin{array}{ccccc}1 & si & x=y=u & ou & x=y=v \\-1 & si & x=u\;et\;y=v & ou & x=v\;et\;y=u \\0 & si & x\notin\{u,v\} & ou & y\notin\{u,v\},\end{array}\right.$$
tandis que si $s(u,v)=1$:
$$(A_e^{(s)})_{xy}=\left\{\begin{array}{ccccc}1 & si & x\in\{u,v\} & et & y\in\{u,v\} \\0  & si & x\notin\{u,v\} & ou & y\notin\{u,v\}.\end{array}\right.$$
Donc
$$d.1_n + A_X^{(s)}=\sum_{e\in E} A_e^{(s)}.$$

On munit l'ensemble des $2^m$ signages de $X$ de la mesure de probabilit\'e uniforme. Pour $e\in E$, on d\'efinit la variable al\'eatoire $A_e:s\mapsto A_e^{(s)}$ \`a valeurs dans les op\'erateurs semi-d\'efinis positifs de rang 1. Les variables al\'eatoires $(A_e)_{e\in E}$ sont ind\'ependantes. Par le th\'eor\`eme \ref{MSS-posproba}, on a pour au moins une r\'ealisation des $A_e$, c-\`a-d. pour au moins un signage $s$:
\begin{equation}\label{signing}
ZM(p_{d.1_n+A_X^{(s)}})\leq ZM(\EE p_{\sum_{e\in E}A_e}).
\end{equation}
Notons $\lambda_{max}^{(s)}$ la plus grande valeur propre de $A_X^{(s)}$. Le membre de gauche de (\ref{signing}) est $d+\lambda_{max}^{(s)}$. Pour le membre de droite, remarquons que 
$$(\EE p_{\sum_{e\in E}A_e})(z)=(\sum_s \frac{1}{2^m}p_{d.1_n+A_X^{(s)}})(z)=(\sum_s \frac{1}{2^m}p_{A_X^{(s)}})(z-d)=(\EE_s p_{A_X^{(s)}})(z-d).$$
Pour continuer, nous aurons besoin de la notion de mariage: un {\it mariage} dans un graphe~$X$ est un ensemble d'ar\^etes 2 \`a 2 disjointes. Notons $p_r$ le nombre de mariages \`a $r$ ar\^etes dans $X$ (on convient que $p_0=1$), et d\'efinissons le {\it polyn\^ome de mariage} de $X$:
$$\mu_X(z)=\sum_{r\geq 0} (-1)^r p_r z^{n-2r}.$$
On a alors un r\'esultat de Godsil et Gutman \cite{GoG} (voir aussi l'appendice A de \cite{MSS-Ra} pour une preuve courte):

\begin{prop}\label{GoGu} $\EE_s p_{A_X^{(s)}}=\mu_X.$
\end{prop}

Comme $\EE p_{\sum_{e\in E} A_e}$ est un polyn\^ome caract\'eristique mixte (voir la proposition \ref{Cor4+15} ci-dessus), on en tire que $\mu_X$ a toute ses racines r\'eelles (ce qui \'etait connu, voir \cite{God}, corollaire 1.2 du chapitre 6). Mieux (m\^eme r\'ef\'erence):
\begin{prop}\label{match} Si $X$ est $d$-r\'egulier, toute racine de $\mu_X$ est r\'eelle et born\'ee par $2\sqrt{d-1}$ en valeur absolue.
\end{prop}
Les propositions \ref{GoGu} et \ref{match} montrent que le membre de droite de (\ref{signing}) satisfait $ZM(\EE p_{\sum_{e\in E}A_e})=d+ZM(\mu_X)\leq d +2\sqrt{d-1}$. L'in\'equation (\ref{signing}) se r\'e\'ecrit donc: 
$$\lambda_{max}^{(s)}\leq 2\sqrt{d-1},$$
et comme le 2-rel\`evement $\tilde{X}^{(s)}$ est bi-colorable, toutes les valeurs propres de $A_X^{(s)}$ sont major\'ees par $2\sqrt{d-1}$ en valeur absolue.
\hfill$\square$

\medskip
Dans \cite{Sr1} on trouvera une application des th\'eor\`emes \ref{MSS-expect} et \ref{MSS-posproba} \`a un probl\`eme de {\oge {\it spectral sparsification}\feg} (\'elagage spectral~?) en th\'eorie des graphes. Ce r\'esultat am\'eliore les r\'esultats ant\'erieurs de Batson-Spielman-Srivastava \cite{BSS} (voir aussi \cite{Naor}).

\bigskip
{\bf Remerciements:} Merci \`a J. Anderson, J. Bourgain, P. Casazza, N. Srivastava et N. Weaver pour d'int\'eressants \'echanges. Je remercie sp\'ecialement Adam Marcus pour son mini-cours sur les entrelacements de polyn\^omes, et Sorin Popa et Stefaan Vaes pour leur relecture attentive d'une premi\`ere version de ce texte. Bravo \`a Terry Tao pour son blog.

\pagebreak
\noindent
Adresse de l'auteur:\\
Universit\'e de Neuch\^atel\\
Facult\'e des Sciences\\
Institut de Math\'ematiques\\
Rue \'Emile--Argand 11\\
CH--2000 Neuch\^atel, Suisse\\
alain.valette@unine.ch

\end{document}